\documentclass[12pt]{amsart}
\newcommand{\CC}{\mathbb{C}}

\newcommand{\NN}{\mathbb{N}}

\newcommand{\RR}{\mathbb{R}}

\newtheorem{definition}{\sc Definition}[section]
\newtheorem{teo}{\sc Theorem}[section]
\newtheorem{prop}{\sc Proposition}[section]

\newtheorem{eje}{\sc Example}[section]
\newtheorem{lemma}{\sc Lemma}[section]
\newtheorem{obs}{\sc Remark}[section]

\begin{document}

\title{Controllability of the Laguerre and the Jacobi Equations.}

\author{DIOMEDES BARCENAS$^{(1)}$, HUGO LEIVA$^{(1)}$, YAMILET QUINTANA$^{(2)}$ \and WILFREDO URBINA$^{(1)}$}
\date{October, 2006.\\
$(1)\,\,\,${\rm Research partially supported by ULA and FONACIT\#G-97000668}\\
$(2)\,\,\,${\rm Research partially supported by  DID-USB under  Grant DI-CB-015-04}}
\email{barcenas@ula.ve, hleiva@ula.ve,yquintana@usb.ve, wurbina@euler.ciens.ucv.ve}
\begin{abstract}
In this paper we study the controllability of  the controlled Laguerre equation and the controlled Jacobi equation. For each case, we found conditions which guarantee when such systems are approximately controllable on the interval $[0, t_1]$. Moreover, we show that these systems can never be exactly controllable.\\

\hspace{-.6cm} {\it Key words and phrases.} Laguerre equation, Jacobi equation, controllability, compact semigroup.\\

\hspace{-.6cm} {\it 2001 Mathematics Subject Classification.}
Primary 93B05. Secondary 93C25.
\end{abstract}

\maketitle
\markboth{\hspace{1cm} DIOMEDES BARCENAS, HUGO LEIVA, YAMILET QUINTANA AND WILFREDO URBINA  \mbox{} }{CONTROLLABILITY OF THE LAGUERRE EQUATION AND THE JACOBI EQUATION}

\section{Introduction.}

The study of orthogonal polynomials which are eigenfunctions of a differential operator
have a long history. In 1929 S. Bochner  \cite{B1} posed the problem of determining all families of  orthogonal polynomials in $\RR$ that are eigenfunctions of some arbitrary but fixed second-order differential operators. In that article, he proved that this property characterizes the so-called classical orthogonal polynomials, linked with the names of Hermite, Laguerre and Jacobi (this last family containing as particular cases the Legendre, Tchebychev and Gegenbauer polynomials). Later 
H.L. Krall and O. Frink \cite{KRO} considered the Bessel polynomials, that are also orthogonal polynomials that satisfies a second order equation, but their orthogonality measure does not have support is $\RR$ but on the unit circle of the complex plane. The general problem, for a differential operator of any order was possed by H. L. Krall  \cite{KR1} in 1938,  he proved that the differential operator has to be of even order and, in \cite{KR2}, he obtained a complete classification for the case of an operator of order four (see \cite{B2}, \cite{KR1},  \cite{KR2} and \cite{mi} for a more detailed references and further developments). There have been recent developments in the direction of connecting the study of orthogonal polynomials with modern problems related to Harmonic Analysis and PDE's, see for instance \cite{BLU},  \cite{GLLNW}, \cite{TOR} .

On the other hand, it is well known that  many differential equations can be solved using the separation variable method, obtaining  solutions in terms of a orthogonal expansion. Nevertheless, is an absolute  merit of C. Sturm and J. Liouville in the 1830s,  the knowledge of the existence of such solutions - long before the advent of Hilbert spaces Theory in the XX -th century-. Their results were precursors of the Operator Theory, but from our present viewpoint can be more  naturally obtained as consequences of the spectral Theorem for compact hermitian operators (the reader is referred to \cite{Y} for the proof of this statement).

With respect to recent developments in controllability of evolution equations of fluid mechanics and controllability of the wave and heat equations via numerical approximation schemes, we refer to \cite{Ima} and \cite{Zua}, respectively.

Following the point of view of connecting the study of diverse aspects of Orthogonal Polynomials Theory with PDE's, in this paper  we are going to study: 

\begin{enumerate}
\item The controllability of controlled Laguerre equation
\begin{equation}
\label{C1}
z_t=\sum_{i=1}^d \Bigg[x_i \frac{\partial^{2} z}{\partial x_i^2} + (\alpha_i+1-x_i) \frac{\partial z}{\partial x_i}\Bigg] +
\sum_{n=0}^{\infty}\sum_{|\nu |= n}u_{\nu}(t)
\langle b,l^{\alpha}_{\nu}\rangle_{\mu_{\alpha}} l^{\alpha}_{\nu}, \ \ t >0, \ \ x \in
\RR_{+}^d,
\end{equation}
 where $\{l^{\alpha}_{\nu}\}$ are the normalized Laguerre polynomials of type $\alpha$ in $d$ variables which are orthogonal polynomials with respect to the the Gamma measure in $\RR_{+}^d$, $\mu_{\alpha}(x) = \prod^{d}_{i=1} \frac{x^{\alpha_i}_{i} e^{-x_i}}{\Gamma(\alpha_i+1)} \; dx$,
$b \in L^{2}(\RR_{+}^d, \mu_{\alpha})$ and the control $u \in L^{2}(0, t_1;
l^{2})$, where  with $l^{2}$ the Hilbert space complex square sumable sequences, that for convenience, it will be written as
 $$l^{2} = \left\{ U = \{
\{U_{\nu} \}_{|\nu |= n} \}_{n \geq 0}: \ \ U_{\nu} \in \CC, \ \
\sum_{n=0}^{\infty}\sum_{|\nu |= n} |U_{\nu} |^{2} < \infty \right\}, $$
with the inner product and norm defined as 

$$ \left\langle U,V \right\rangle_{l^2} = \sum_{n=0}^{\infty}\sum_{|\nu|=
n}U_{\nu}\overline{V_{\nu}}, \, \ \ \|U \|_{l^{2}}^2 =
\sum_{n=0}^{\infty}\sum_{|\nu |=n}|U_{\nu} |^2, \ \ U, V \in
l^{2}. $$ 

We will prove the following statement: If for all $\nu =
(\nu_1, \nu_2, \ldots, \nu_d)\in \NN_{0}^d$ $$
\langle b,l^{\alpha}_{\nu}\rangle_{\mu_{\alpha}} = \int_{\RR_{+}^d}b(x)
l^{\alpha}_{\nu}(x)\mu_{\alpha}(dx) \not=0, $$  then the system is
approximately controllable on $[0, t_1].$ Moreover, the system can
never be exactly controllable.

In particular, we consider the Laguerre
equation in one variable with a single control $$ z_{t}  =  x z_{xx} +(\alpha +1-x) z_x + b(x)u \
\ t \geq 0, \ \  x \in \RR_{+}, $$ where $b \in L^{2}(\RR_{+},\mu_{\alpha})$ and
the control $u$ belong to $L^{2}(0,t_1;\RR_{+})$. This system is
approximately controllable if and only if
$$ \int_{\RR_{+}} b(x) l^{\alpha}_{\nu}(x)x^{-\alpha}e^{x} dx \not= 0, \ \ \nu = 0,1,2,\dots. $$

\item The controllability of controlled Jacobi equation
\begin{equation}
\label{C2}
z_t=\sum_{i=1}^d \Bigg[(1-x^2_i) \frac{\partial^{2} z}{\partial x_i^2} + (\beta_{i} -\alpha_{i}-\left(\alpha_{i} +\beta_{i} +2\right)x_i) \frac{\partial z}{\partial x_i}\Bigg] +
\sum_{n=0}^{\infty}\sum_{|\nu |= n}u_{\nu}(t)
\langle b,p^{\alpha,\beta}_{\nu}\rangle_{\mu_{\alpha,\beta}} p^{\alpha,\beta}_{\nu}, 
\end{equation}
$ t >0, \, x \in [-1,1]^d$
 where $\{p^{\alpha,\beta}_{\nu}\}$ are the normalized Jacobi polynomials of type $\alpha = (\alpha_1,\ldots,\alpha_d),\, \beta = (\beta_1,\ldots,\beta_d)\in \RR^d$, $\alpha_i,\beta_i> -1$, in $d$ variables, which are orthogonal polynomials with respect to the Jacobi measure in $[-1,1]^d$$\mu_{\alpha,\beta}(x) = \prod^{d}_{i=1}(1-x_i)^{\alpha_i}(1+x_i)^{\beta_i} \; dx$, $b \in L^{2}([-1,1]^d, \mu_{\alpha,\beta})$ and the control $u \in L^{2}(0, t_1;
l^{2})$.

Analogous to the previous case, we will prove that if for all $\nu =
(\nu_1, \nu_2, \dots, \nu_d)\in \NN_{0}^d$ $$
\langle b,p^{\alpha,\beta}_{\nu}\rangle_{\mu_{\alpha,\beta}} = \int_{[-1,1]^d}b(x)
p^{\alpha,\beta}_{\nu}(x)\mu_{\alpha,\beta}(dx) \not=0, $$  then the system is
approximately controllable on $[0, t_1]$; but, it  can never be exactly controllable.

Also, in particular, for  $\alpha ,\beta >-1$ we consider the Jacobi
equation  in one variable with a single control $$ z_{t}  =  (1-x^2) z_{xx} +(\beta -\alpha-\left(\alpha +\beta +2\right)x) z_x + b(x)u,
\, t \geq 0, \ \  x \in [-1,1], $$ where $b \in L^{2}([-1,1],\mu_{\alpha,\beta})$ and
the control $u$ belong to $L^{2}(0,t_1;[-1,1])$. This system is
approximately controllable if and only if 
$$  \int_{[-1,1]} b(x)p^{\alpha,\beta}_{\nu} \left( 1-x\right) ^{-\alpha }\left( 1+x\right) ^{-\beta
}dx \not= 0, \,  \nu = 0,1,2,3,\dots.$$
\end{enumerate}

The Laguerre differential operator,
\begin{eqnarray}
\mathcal{L}^{\alpha} &=& - \sum^{d}_{i=1} \Bigg[ x_i \partial^2_{x_i}
+ (\alpha_i + 1 - x_i ) \partial_{x_i} \Bigg]\\
\nonumber \mbox{and the Jacobi differential operator,}\\
\mathcal{L}^{\alpha,\beta} &=& - \sum^{d}_{i=1} \Bigg[ (1-x_i^2)\partial^2_{x_i}
+  (\beta_{i} -\alpha_{i}-\left(\alpha_{i} +\beta_{i} +2\right)x_i) \partial_{x_i} \Bigg]
\end{eqnarray}
are well-known operators in the theory Orthogonal Polynomials , in Probability Theory, 
in Quantum Mechanics  and in Differential Geometry (see \cite{GLLNW}, \cite{me1}, \cite{mi}, \cite{mu1},\cite{Sz}). 

With the results of this paper, we complete the study of  controllability problem for the operators associated to classical orthogonal polynomials. In a previous paper \cite{BLU} it was considered  the case of Ornstein-Uhlenbeck operator, and as far as we know, these controled equations have not been studied until now. Also we obtain  results, as in \cite{DR}, on approximate controllability for some higher dimensional systems associated to a Sturm-Liouville operators of the form
\[
\mathcal{L} =\frac{1}{\rho(x)}\sum_{i,j=1}^{d}\partial_{ x^{i}}\left( a_{j}^{i}(x)\partial_{x^{j}}\right),
\]
 where $x\in \RR^{d}$, $\rho:\RR^{d}\rightarrow \RR$ is a constant function and $A(x)=\left( a_{j}^{i}(x)\right)_{1\leq i,j\leq d}$ is a constant matrix. It remains open the study of the general case.  The  arguments used in this paper can be extended to this more general setting.

Two important tools which allow to improve and complete the study of controllability problem for the operator associated to classical orthogonal polynomials were used in \cite{BLU} and come from \cite{BLS} (Theorem 3.3) and 
\cite{CP1} (Theorem  A.3.22).

The outline of the paper is the following. Section 2 is dedicated to preliminary results. Section 3 we present  main results of the paper, the controllability of the controlled Laguerre equation (\ref{C1}) and the  controllability of the controlled Jacobi equation (\ref{C2}).

\section{Preliminary  results.}

In this section we shall choose the spaces where our problems will
be set and we shall present some results that are needed  in the next section. Also,we will give the definition of exact and approximate controllability. 

To deal with polynomials in several variables we use the standard multi-index notation. A multi-index is denoted by $\nu=(\nu_1,\ldots,\nu_d)\in \NN_{0}^d$, where $\NN_{0}$ is the set of non negative integers numbers. For $\nu\in\NN_{0}^{d}$ we denote by $\nu
!=\prod_{i=1}^d\nu _i!$,  $\left| \nu \right| =\sum_{i=1}^d\nu _i$,  $
\partial _i=\frac \partial {\partial x_i},$ for each $1\leq i\leq d$ and $
\partial ^\nu =\partial _1^{\nu _1}\ldots \partial _d^{\nu _d}$.

Then  the normalized Laguerre polynomials of type $\alpha = (\alpha_1 ,\ldots,\alpha_d)\in \RR^d$, $\alpha_i> -1$, and order $\nu$ in $d$ variables is given by the tensor product
\begin{equation}
l^{\alpha}_{\nu} (x)={\sqrt{\nu!} \over
\sqrt{\Gamma(\alpha+\nu+1)}}\prod_{i=1}^d(-1)^{\alpha _i}x_i^{-\alpha_{i}}e^{x_i}\frac{\partial ^{\nu_i}}{\partial x_i^{\nu _i}}(x_i^{\nu_i+\alpha_{i}} e^{-x_i}).
\end{equation}
It is well known, that the Laguerre polynomials are eigenfunctions
of the Laguerre operator $\mathcal{L}^{\alpha}$,
\[
\mathcal{L}^{\alpha} l^{\alpha}_{\nu}(x)=-\left|\nu \right| l^{\alpha}_{\nu}(x).
\]
Given a function $f \in L^{2}(\RR_{+}^d, \mu_{\alpha})$ its
$\nu$-Fourier-Laguerre coefficient is defined by
\[
\langle f,l^{\alpha}_{\nu} \rangle_{\mu_{\alpha}} =\int_{\RR_{+}^d} f(x) l^{\alpha}_{\nu} (x)\mu_{\alpha}(dx),
\]
Let $C^{\alpha}_n$ be the closed subspace of $L^{2}(\RR_{+}^d, \mu_{\alpha})$ generated by $\left\{ l^{\alpha}_{\nu} \ :\left| \nu
\right| =n\right\}$, $C^{\alpha}_n$ is a finite dimensional subspace of
dimension ${n+d-1 \choose n}$. By the ortogonality of the Laguerre
polynomials with respect to $\mu_{\alpha}$ it is easy to see that
$\{C^{\alpha}_n\}$ is a orthogonal decomposition of $L^{2}(\RR_{+}^d, \mu_{\alpha})$, $$
L^{2}(\RR_{+}^d, \mu_{\alpha}) = \bigoplus_{n=0}^{\infty} C^{\alpha}_n ,$$ which is called
the Wiener-Laguerre chaos.

The orthogonal projection $P^{\alpha}_n$ of $L^{2}(\RR_{+}^d, \mu_{\alpha})$ onto $C^{\alpha}_n$ is given by 
$$ P^{\alpha}_n f=\sum_{\left|\alpha\right|=n}\langle f,l^{\alpha}_{\nu}\rangle_{\mu_{\alpha}} l^{\alpha}_{\nu},\ \ f \in L^{2}(\RR_{+}^d, \mu_{\alpha}),
$$
and for a given $f\in L^{2}(\RR_{+}^d, \mu_{\alpha})$ its Laguerre expansion is
given by $f= \sum_n P^{\alpha}_n f.$

Using this notation one can prove the following  espectral
decomposition of $\mathcal{L}^{\alpha}$ 
$$\mathcal{L}^{\alpha}f = \sum_{n=0}^{\infty} (-n) P^{\alpha}_nf, \ \ f \in
L^{2}(\RR_{+}^d, \mu_{\alpha}), $$ and its domain $D(\mathcal{L}^{\alpha})$ is 
$$ D(\mathcal{L}^{\alpha})=\left\{f\in
L^{2}(\RR_{+}^d, \mu_{\alpha}): \sum_{n=0}^{\infty} n^2 \|P^{\alpha}_nf\|_{2, \mu_{\alpha}}< \infty\right\}.$$

Let $Z = L^{2}(\RR_{+}^d, \mu_{\alpha})$ and $l^{2}$ be the Hilbert space of  complex square sumable sequences. Now, suppose that $b$ is a fixed element of
$Z$ and consider the linear and bounded operator $B
:l^{2}\rightarrow Z$ defined by
\begin{equation}
\label{F1}
BU =\sum_{n=0}^{\infty}\sum_{|\nu |= n}U_{\nu}
\langle b,l^{\alpha}_{\nu} \rangle_{\mu_{\alpha}}l^{\alpha}_{\nu}.
\end{equation}
Then, the system (\ref{C1}) can be written as follows
\begin{equation}\label{C11}
z'= \mathcal{L}^{\alpha} z + Bu, \ \ t >0.
\end{equation}

By a similar way, the normalized Jacobi polynomials of type $\alpha = (\alpha_1,\ldots,\alpha_d), \beta = (\beta_1,\ldots,\beta_d)\in \RR^d$, $\alpha_i,\beta_i> -1$, of order $\nu$ in $d$ variables is given by the tensor product
\begin{equation}
p^{\alpha,\beta}_{\nu} (x)=(h_{\nu}^{\left( \alpha ,\beta \right)})^{-1/2}\prod_{i=1}^d\left( 1-x_i\right) ^{-\alpha_i }\left( 1+x_i\right) ^{-\beta_i
} \frac{\left( -1\right) ^{\nu_i}}{2^{\nu_i}\nu_i!}\frac{d^{\nu_i}}{dx_i^{\nu_i}}\left\{
\left( 1-x_i\right) ^{\alpha_i +\nu_i}\left( 1+x_i\right) ^{\beta_i +\nu_i}\right\},
\end{equation}

where $ h_{\nu}^{\left( \alpha ,\beta \right)} = \prod_{i=1}^d  h_{\nu_i}^{\left( \alpha_i ,\beta_i \right)},$
with $h_{\nu_i}^{\left( \alpha_i ,\beta_i \right)}  = \frac{2^{\alpha_i +\beta_i +1}}{2\nu_i+\alpha_i+\beta_i +1} \frac{\Gamma \left( \nu_i+\alpha_i +1 \right) \Gamma \left( \nu_i+\beta_i +1\right)}{\Gamma \left( \nu_i+1\right) \Gamma \left( \nu_i+\alpha_i +\beta_i +1 \right) }$.

Also, it is well-known, that the Jacobi polynomials are eigenfunctions
of the Jacobi operator $\mathcal{L}^{\alpha,\beta}$,
\[
\mathcal{L}^{\alpha,\beta} p^{\alpha,\beta}_{\nu} = - \sum^{d}_{i=1} \bigg[ (1-x_i^2) \partial^2_{x_i} p^{\alpha,\beta}_{\nu}
+  (\beta -\alpha-\left(\alpha +\beta +2\right)x_i) \partial_{x_i}p^{\alpha,\beta}_{\nu} \bigg] 
= \sum_{i=1}^d \nu_i\left( \nu_i+\alpha_i +\beta_i +1\right) p^{\alpha,\beta}_{\nu}. 
\]
And given a function $f \in L^{2}([-1,1]^d, \mu_{\alpha,\beta})$ its
$\nu$-Fourier-Jacobi coefficient is defined by
\[
\left\langle f,p^{\alpha,\beta}_{\nu} \right\rangle_{\mu_{\alpha,\beta}} =\int_{[-1,1]^d} f(x) p^{\alpha,\beta}_{\nu} (x)\mu_{\alpha,\beta}(dx).
\]

As the eigenvalues of the Jacobi
operator are not linear in $n$, following \cite{BU} we are going to consider a alternative decomposition, in order to obtain an espectral decomposition  of $\mathcal{L}^{\alpha,\beta}  f$ for any $f \in L^2([-1,1]^d, \mu_{\alpha, \beta} )$ in
terms of the orthogonal projections.

For fixed $\alpha = ( \alpha_1, \alpha_2, \cdots, \alpha_d),$
$\beta = ( \beta_1, \beta_2, \cdots, \beta_d)$,  in $\mathbb{R}^d$
such that $\alpha_i, \beta_i > {-{1\over 2}}$ let us consider the set,
$$R^{\alpha, \beta} =\left\{ r\in \mathbb{R}^+
: \mbox{there exists}\,(\kappa_1,\ldots,\kappa_n) \in \mathbb{N}^d_0,  \mbox{with} \,
r=\sum_{i=1}^d \kappa_i(\kappa_i+\alpha_i+\beta_i +1) \right\}.$$
$R^{\alpha, \beta} $ is a numerable subset of $\mathbb{R}^+$, we  can write an
enumeration of  $R^{\alpha, \beta} $ as 
$\{r_n\}_{n=0}^\infty$ with $0= r_0<r_1<\cdots. $ Let
$$A_n^{\alpha, \beta} =\left\{ \kappa=(\kappa_1,\ldots,\kappa_d) \in
\mathbb{N}^d_0: \sum_{i=1}^d \kappa_i(\kappa_i+\alpha_i+\beta_i+1) =
r_n \right\}.$$ 
Notice that $A_0^{\alpha, \beta} =\{(0,\ldots,0)\}$ and that if $\kappa\in A_n^{\alpha, \beta} $ then $\sum_{i=1}^d
\kappa_i(\kappa_i+\alpha_i+\beta_i+1) =r_n.$

Let $C_n^{\alpha, \beta} $ denote the closed subspace of
$L^2([-1,1]^d,\mu_{\alpha, \beta} )$ generated by the linear combinations of
$\{ p_\kappa^{\,\,\alpha, \beta}: \kappa\in A_n^{\alpha, \beta}  \}.$ By the orthogonality of the
Jacobi polynomials with respect to $\mu_{\alpha, \beta}$ and the density of the polynomials, it is not difficult
to see that $\{C^{\alpha, \beta} _n\}$ is an orthogonal decomposition of
$L^2([-1,1]^d,\mu_{\alpha, \beta})$, that is
\begin{equation}\label{Caosdecomp2}
 L^2([-1,1]^d,\mu_{\alpha, \beta})=\bigoplus_{n=0}^\infty C^{\alpha, \beta} _n. 
\end{equation}

We call (\ref{Caosdecomp2}) a modified Wiener--Jacobi
decomposition.

The ortogonal proyection $ P^{\alpha,\beta}_n$ of $L^2([-1,1]^d,\mu_{\alpha, \beta})$ onto $C_n^{\alpha, \beta} $ is given by 
$$ P^{\alpha,\beta}_n f=\sum_{\nu \in A_n^{\alpha,\beta}}\langle f,p^{\alpha,\beta}_{\nu}\rangle_{\mu_{\alpha,\beta}} p^{\alpha,\beta}_{\nu},\ \ f \in L^{2}([-1,1]^d, \mu_{\alpha,\beta}),
$$
and for a given $f\in L^{2}([-1,1]^d, \mu_{\alpha,\beta})$ its Jacobi expansion is
then given by  $$f= \sum_{n=0}^{\infty} P^{\alpha,\beta}_n f.$$
Therefore  $\left\{ P^{\alpha,\beta}_n\right\} _{n \geq 0}$ is a complete system of orthogonal projections in $L^{2}([-1,1]^d, \mu_{\alpha,\beta})$.

Using this notation one can prove the following  espectral
decomposition of the operator $\mathcal{L}^{\alpha,\beta}$ 
$$\mathcal{L}^{\alpha,\beta} = \sum_{n=0}^{\infty} (-r_n) P^{\alpha,\beta}_n f ,$$
 $ f \in L^{2}([-1,1]^d, \mu_{\alpha,\beta}), $ and its domain $D(\mathcal{L}^{\alpha,\beta})$ is given by
$$ D(\mathcal{L}^{\alpha,\beta})=\left\{ f\in L^{2}([-1,1]^d, \mu_{\alpha,\beta}):  \sum_{n=0}^{\infty}
(r_n)^2 \|P^{\alpha,\beta}_nf\|_{2, \mu_{\alpha, \beta}}< \infty\right\}.$$\\

Let $W = L^{2}([-1,1]^d, \mu_{\alpha,\beta})$ and $l^{2}$ be the Hilbert space of  complex square sumable sequences. Again, suppose that $b$ is a fixed element of $W$ and consider the linear and bounded operator $\tilde{B}:l^{2}\rightarrow W$ defined by
\begin{equation}
\label{j1}
\tilde{B}U =\sum_{n=0}^{\infty}\sum_{|\nu |= n}U_{\nu}
\langle b,p^{\alpha,\beta}_{\nu}\rangle_{\mu_{\alpha,\beta}}p^{\alpha,\beta}_{\nu}.
\end{equation}
Then, the system (\ref{C2}) can be written as follows
\begin{equation}
\label{C21}
w'= \mathcal{L}^{\alpha,\beta} w + \tilde{B}\tilde{u}, \ \ t >0,
\end{equation}

\begin{teo}
\label{T1}
The operators $\mathcal{L}^{\alpha}$ and  $\mathcal{L}^{\alpha,\beta}$ are  the infinitesimal generators of  analytic
semigroups $\left\{ T^{\alpha}(t) \right\}_{t \geq 0}$ and $\left\{ T^{\alpha,\beta}(t) \right\}_{t \geq 0}$, respectively. They are given as 
\begin{equation}
\label{C4}
T^{\alpha}(t)z=\sum_{n=0}^{\infty}e^{-nt}P^{\alpha}_n z, \ \ z \in Z, \ \ t \geq 0,
\end{equation}
where  $\left\{ P^{\alpha}_n\right\} _{n \geq 0}$ is a complete orthogonal
projections in the Hilbert space $Z$ given by\\
$P^{\alpha}_n z=\sum_{\left|\nu\right|=n}\left \langle z,l^{\alpha}_{\nu}
\right\rangle_{\mu_{\alpha}}l^{\alpha}_{\nu}, \ \ n \geq 0, \ \ z \in Z,$ and
\begin{equation}
\label{C41}
T^{\alpha,\beta}(t)w=\sum_{n=0}^{\infty}e^{-r_nt}P^{\alpha,\beta}_n w, \ \ w \in W, \ \ t \geq 0,
\end{equation}
where  $\left\{ P^{\alpha,\beta}_n\right\} _{n \geq 0}$ is a complete orthogonal
projections in the Hilbert space $W$ given by\\
$P^{\alpha,\beta}_n w=\sum_{\nu \in A_n^{\alpha,\beta}}\langle w,p^{\alpha,\beta}_{\nu}\rangle_{\mu_{\alpha,\beta}} p^{\alpha,\beta}_{\nu}, \ \ n \geq 0, \ \ w \in W.$
\end{teo}

\begin{lemma}
\label{comp} The semigroups given by (\ref{C4}) and (\ref{C41})
are compact for $t>0$.
\end{lemma}

\begin{proof}

Since $T^{\alpha}(t)$ is given by $$ T^{\alpha}(t)z = \sum_{n=0}^{\infty}
e^{-nt}P^{\alpha}_n z, \ \ \ \ t >0, $$ 
we can consider the following
sequence of compact operators
$$ T^{\alpha}_{k}(t)z = \sum_{n=0}^{k} e^{-n
t}P^{\alpha}_n z, \ \ \ \ t > 0. $$ 

It is easy to see that the sequence of compact operators $\{T^{\alpha}_{n}(t)\}$ converges uniformly to $T^{\alpha}(t)$ for all $t >0$. 

Analogously, $T^{\alpha,\beta}(t)$ is given by 
$$ T^{\alpha,\beta}(t)w = \sum_{n=0}^{\infty}e^{-r_nt}P^{\alpha,\beta}_n w, \ \ \ \ t >0, $$
so that, we can consider the following
sequence of compact operators 
$$ T^{\alpha,\beta}_{k}(t)w =\sum_{n=0}^{k}e^{-r_nt}P^{\alpha,\beta}_n w, \ \ \ \ t > 0. $$ 
and again it is easy to see that the sequence of compact operators $\{T^{\alpha,\beta}_{k}(t)\}$ converges uniformly to $T^{\alpha,\beta}(t)$ for all $t >0$.

Then, from part e) of Theorem A.3.22 of \cite{CP1} we conclude the compactness of the semigroups 
$T^{\alpha}(t)$ and $T^{\alpha,\beta}(t)$, respectively.
\end{proof}

Now, we shall give the definitions of exact and approximate
controllability in terms of system (\ref{C11}) and (\ref{C21}).  In spite of this definitions can be given for more general evolutions equations, we concentrated our atention to the cases of our interest. 

For all $z_0 \in Z$, $w_0\in W$ and  given controls $u \in L^{2}(0,t_1;
l^{2})$ and $\tilde{u}\in  L^{2}(0,t_1;
l^{2})$ the equations (\ref{C11}) and (\ref{C21})  have a unique mild solution given -in each case- by
\begin{equation}
\label{C12}
z(t) = T^{\alpha}(t)z_0 + \int_{0}^{t} T^{\alpha}(t-s)Bu(s)ds, \ \ 0 \leq t \leq
t_1.
\end{equation}
\begin{equation}
\label{C22}
w(t) = T^{\alpha,\beta}(t)w_0 + \int_{0}^{t} T^{\alpha,\beta}(t-s)\tilde{B}\tilde{u}(s)ds, \ \ 0 \leq t \leq
t_1.
\end{equation}

\begin{definition} 
\label{D:d1}
(Exact Controllability).\\
 We shall say that the system
(\ref{C11}) (respectively, (\ref{C21})) is exactly controllable on $[0,t_1], \ \ t_1>0$, if
for all $z_0, z_1 \in Z$ (respectively, $w_0, w_1 \in W$)
 there exists a control $u \in L^{2}(0,t_1;l^{2})$\\(respectively, $\tilde{u}\in  L^{2}(0,t_1;l^{2})$ )
 such that the solution $z(t)$ of (\ref{C11}) corresponding to $u$\\(respectively, the solution $w(t)$ of (\ref{C21}) corresponding to $\tilde{u}$), that verifies $z(t_1) = z_1$ (respectively, $w(t_1 )=w_1$).
\end{definition}

Consider the following bounded linear operators
\begin{equation} 
\label{c5}
G : L^{2}(0,t_1;l^{2}) \rightarrow Z, \ \ Gu
=\int_{0}^{t_1} T^{\alpha}(t_1-s) Bu(s) ds,
\end{equation}
\begin{equation} 
\label{c51}
\tilde{G} : L^{2}(0,t_1;l^{2}) \rightarrow W, \ \ \tilde{G}\tilde{u}
=\int_{0}^{t_1} T^{\alpha,\beta}(t_1-s) \tilde{B}\tilde{u}(s) ds.
\end{equation}

Then, the following Proposition is a characterization of the exact
controllability of the sytems (\ref{C11}) and (\ref{C21}).
\begin{prop} 
\hspace{1cm}
\label{p1}
\begin{enumerate}
\item[i)] The system (\ref{C11}) is exactly controllable on $[0,t_1]$ if and
only if, the operator $G$ is surjective, that is to say $$ G
L^{2}(0,t_1;l^{2}) = G L^{2}= \mbox{Range}(G)= Z. $$
\item[ii)]The system (\ref{C21}) is exactly controllable on $[0,t_1]$ if and
only if, the operator $\tilde{G}$ is surjective, that is to say 
$$ \tilde{G}L^{2}(0,t_1;l^{2}) = \tilde{G} L^{2}= \mbox{Range}(\tilde{G})= W. $$
\end{enumerate}
\end{prop}

\begin{definition}
\label{DD1}
We say that (\ref{C11})\, (respectively, (\ref{C21})) is approximately controllable in $[0,t_1]$
if for all $z_0, z_1 \in Z$ (respectively, $w_0, w_1 \in W$) and $\epsilon >0$, there exists a
control $u \in L^{2}(0,t_1; l^{2})$ (respectively, $\tilde{u}\in  L^{2}(0,t_1;l^{2})$) such that the
solution $z(t)$  given by (\ref{C12}) (respectively, the solution $w(t)$ given by (\ref{C22})) satisfies 
$$ \|z(t_1) -z_1 \|\leq \epsilon, \,\, (\mbox{respectively, }\|w(t_1) -w_1 \|\leq \epsilon). $$
\end{definition}

 Via duality, the following Theorem allows to  give a characterization of the approximate controllability for our systems. Such characterization  holds in general and  the reader is referred to \cite{CP1} for the details of its  proof.

\begin{teo} 
\label{T2}
\hspace{1cm}
\begin{enumerate}
\item[i)]
The system (\ref{C11})  is approximately controllable on $[0,t_1]$ if and only if
\begin{equation}\label{C13}
B^{*}\left(T^{\alpha}\right)^{*}(t)z =0, \ \ \forall t \in [0,t_1], \mbox{implies} \,  \, z=0.
\end{equation}
\item[ii)] The system (\ref{C21}) is approximately controllable on $[0,t_1]$ if and only if
\begin{equation}
\label{C33}
\tilde{B}^{*}\left(T^{\alpha,\beta}\right)^{*}(t)w =0, \ \ \forall t \in [0,t_1], \mbox{implies} \, \, w=0.
\end{equation}
\end{enumerate}
\end{teo}

 \section{Controllability of the controlled Laguerre equation and the controlled Jacobi equation.}

In this section we shall prove  the main results of the paper, 

\begin{teo}
\label{T3}
\hspace{1cm}
\begin{enumerate}
\item[i)] If for all $n \in \NN_{0}$ and $|\nu|= n$
we have
\begin{equation}
\label{mila}
 \left\langle b,l^{\alpha}_{\nu}\right\rangle_{\mu_{\alpha}} = \int_{\RR_{+}^d}b(x)
l^{\alpha}_{\nu}(x)\mu_{\alpha}(dx) \not=0,
\end{equation}
then the system
(\ref{C11}) is approximately controllable on $[0, t_1]$, but never
exactly controllable.
\item[ii)] If for all $n \in \NN_{0}$ and $|\nu|= n$
we have
\begin{equation}
\label{mila1}
 \langle b,p^{\alpha,\beta}_{\nu}\rangle_{\mu_{\alpha,\beta}} = \int_{[-1,1]^d}b(x)
p^{\alpha,\beta}_{\nu}(x)\mu_{\alpha,\beta}(dx) \not=0,
\end{equation}then the system
(\ref{C21}) is approximately controllable on $[0, t_1]$, but never
exactly controllable.
\end{enumerate}
\end{teo}

\begin{obs}
\label{obs1}
Notice that it is sufficient  to prove the first part of the Theorem, since the proof depends of relation between the adjoint operator of $B$ (respectively, $\tilde{B}$) and the adjoint operator of  $T^{\alpha}(t)$ (respectively, $T^{\alpha,\beta}(t)$) given by the Theorem \ref{T2}.
\end{obs}

\begin{proof} 
Suppose condition  (\ref{mila}). Next, we compute $B^{*}: Z
\rightarrow l^{2}$. In fact,
\begin{eqnarray*}
\left\langle BU,z \right\rangle_{\mu_{\alpha}} & = &
\left\langle\sum_{n=0}^{\infty}\sum_{|\nu |=
n}U_{\nu}\left\langle b,l^{\alpha}_{\nu}\right\rangle_{\mu_{\alpha}} l^{\alpha}_{\nu},z
\right\rangle _{Z,Z}\\ & = & \sum_{n=0}^{\infty}\sum_{|\nu |=
n}U_{\nu}\left\langle b,l^{\alpha}_{\nu}\right\rangle_{\mu_{\alpha}}\left\langle z , l^{\alpha}_{\nu}
\right\rangle_{Z,Z}\\ & = & \left\langle  U, \{ \{
\left\langle b,l^{\alpha}_{\nu}\right\rangle_{\mu_{\alpha}}\left\langle z , l^{\alpha}_{\nu}
\right\rangle\}_{|\nu|= n} \}_{n \geq
0}\right\rangle_{l^{2},l^{2}}.
\end{eqnarray*}
Therefore, 
$$ B^{*}z =  \{ \{
\langle b,l^{\alpha}_{\nu}\rangle_{\mu_{\alpha}}\langle z , l^{\alpha}_{\nu} \rangle
\}_{|\nu|= n} \}_{n \geq 0} = \sum_{n=0}^{\infty}\sum_{|\nu |=
n}\langle b,l^{\alpha}_{\nu}\rangle_{\mu_{\alpha}}\langle z , l^{\alpha}_{\nu}
\rangle e_{\nu}, $$ where $\{ \{ e_{\nu}\}_{|\nu|= n} \}_{n
\geq 0}$ is the canonical basis of $l^{2}$.\\

On the other hand,
$$ (T^{\alpha})^{*}(t)z=\sum_{n=0}^{\infty}e^{-nt}P^{\alpha}_nz, \
\ z \in Z, \ \ t \geq 0. $$
Then,

$$ B^{*}(T^{\alpha})^{*}(t)z = \{ \{
\langle b,l^{\alpha}_{\nu}\rangle_{\mu_{\alpha}}\langle(T^{\alpha})^{*}(t)z , l^{\alpha}_{\nu} \rangle
\}_{|\nu|= n} \}_{n \geq 0}. $$ 

According with the part $i)$ of Theorem \ref{T2}
the system (\ref{C11}) is approximately controllable on $[0,t_1]$
if and only if 
\begin{equation}
\label{C15}
\langle b,l^{\alpha}_{\nu}\rangle_{\mu_{\alpha}}\langle(T^{\alpha})^{*}(t)z , l^{\alpha}_{\nu} \rangle=
0, \ \ \forall t \in [0,t_1], \ \ |\nu |=n, \ n=0,2, \cdots,
\infty, \ \ \Rightarrow z=0.
\end{equation}

Since $\langle b,l^{\alpha}_{\nu}\rangle_{\mu_{\alpha}} \not=0\,$ for $|\nu |=n,\, n
\geq 0$, then condition (\ref{C15}) is equivalent to
\begin{equation}\label{C52}
\langle(T^{\alpha})^{*}(t)z , l^{\alpha}_{\nu} \rangle= 0, \ \ \forall t \in [0,t_1], \
\ |\nu |=n, \, n \geq 0, \ \ \Rightarrow z=0.
\end{equation}

Now, we shall check condition (\ref{C52}):

 $$\langle(T^{\alpha})^{*}(t)z , l^{\alpha}_{\nu} \rangle= \sum_{m=0}^{\infty}e^{-mt}
\langle P_{m}z,l^{\alpha}_{\nu}l^{\alpha}_{\nu} \rangle = 0, \ \ |\nu |=n, \ \ n=0,1, 2, \ldots, \infty;
\ \ t \in [0,t_1]. $$

Applying  Lemma 3.14 from \cite{CP1}, pag. 62 (see also Lemma 3.1 of \cite {BLU}), we conclude that
$$ \langle P^{\alpha}_{m}z,l^{\alpha}_{\nu}\rangle= 0, \ \ |\nu |=n, \ \ m, n=0,1,2, \ldots,
\infty.
 $$
 i.e.,
 $$
\sum_{|\nu |= m}\left\langle z,l^{\alpha}_{\nu}\right\rangle\left\langle l^{\alpha}_{\nu} ,
l^{\alpha}_{\nu} \right\rangle = 0, \ \ |\nu |=n, \ \ m, n=0,1,2, \ldots, \infty.
 $$
 i.e.,
 $$
\left\langle z,l^{\alpha}_{\nu}\right\rangle = 0, \ \ |\nu |=n, \ \ n=0,1,2, \ldots,
\infty.$$

Since $\{  l^{\alpha}_{\nu} \}_{\nu} $ is a complete orthonormal basis of $Z$, we conclude that $z=0$.

On the other hand, from Lemma \ref{comp} we know that $T^{\alpha}(t)$ is compact for $t > 0$,
then applying Theorem 3.3 from \cite{BLS} we conclude that the
system (\ref{C11}) is not exactly controllable on any interval $[0,
t_1]$. This last fact and the remark  \ref{obs1} finish the proof.
\end{proof}

Since an  important ingredient in the above proof is Theorem 3.3 from \cite{BLS},  for completeness of this work we shall include here its  proof -adapted to our context-. 

In fact, from Proposition \ref{p1} it is enough to prove that the operator
$$ G : L^{2}(0,t_1; l^{2}) \rightarrow Z, \ \ Gu
=\int_{0}^{t_1} T^{\alpha}(t_1-s) Bu(s) ds $$ satisfies $$
\mbox{Range}(G)\not= Z. $$ 
In order to do that, we shall prove that the
operator $G$ is compact. For all $\delta >0$ small enough
the operator $G$ can be written as follows 
$$ G = G_{\delta} +
S_{\delta}, \ \ G_{\delta}, S_{\delta} \in L(L^{2}(0,t_1;l^{2}, Z), $$
 where 
 $$ G_{\delta}u = \int_{0}^{t_1 - \delta } T^{\alpha}(t_1-s) Bu(s) ds
 \ \ \mbox{and } \ \ S_{\delta}u=
 \int_{t_1 - \delta }^{t_1} T^{\alpha}(t_1-s) Bu(s) ds.$$

{\bf Claim 1.} The operator $G_{\delta}$ is compact. In fact,
 \begin{eqnarray*}
 G_{\delta}u & = & \int_{0}^{t_1 - \delta }
 T^{\alpha}( \delta )T^{\alpha}(t_1 - \delta -s) Bu(s) ds\\
 & = &  T^{\alpha}(\delta )\int_{0}^{t_1 - \delta }
T^{\alpha}(t_1 - \delta -s) Bu(s) ds\\
 & = & T^{\alpha}(\delta )H_{\delta}u.
 \end{eqnarray*}
 Since $T^{\alpha}(\delta )$ is compact and $H_{\delta} \in
 L(L^{2}(0,t_1; l^{2}), Z)$, then $ G_{\delta}$ is compact.\\

 {\bf Claim 2.} For $\epsilon >0$ there exists $\delta >0$ such that
 $\|S_{\delta} \| < \epsilon$. In fact,
 \begin{eqnarray*}
 \|S_{\delta}u \| & \leq  & \int_{t_1 - \delta}^{t_1  }
 \| T^{\alpha}(t_1 -s)\| \| B \|\|u(s)\| ds\\
 & \leq & \int_{t_1 - \delta}^{t_1  }
 M \| B \|\|u(s)\| ds,
 \end{eqnarray*}
 where
 $$
 M = \sup_{0 \leq s \leq t \leq t_1}\|T^{\alpha}(t-s) \|.
 $$
 Applying H\"older's inequality we obtain
 $$
 \|S_{\delta}u \| \leq M \|B \|\delta \|u \|_{L^{2}}.
 $$
 Therefore,  $\|S_{\delta} \|< \epsilon$ if $\delta <
 {\epsilon \over M \|B \|}$.

 Hence, for all natural  number $n$ the exists $\delta_n >0$ such that
 $$
 \|G -G_{\delta_n} \| = \|S_{\delta_n} \|< {1 \over n }, \ \
 n=1,2,3,\dots.
 $$
So that, the sequence of compact operators $\{G_{\delta_n} \}$
converges uniformly to  $G$. Then applying part e) of Theorem
A.3.22 from \cite{CP1} we obtain that $G$ is compact. Finally, from part g)
of the same Theorem we obtain that $\mbox{Range}(G)\not= Z$.

As special cases of Theorem \ref{T3} we consider
\begin{eje}
\hspace{1cm}
\begin{enumerate}
\item[a)] The Laguerre equation in one variable with a single control
\begin{equation}
\label{C111}
z_{t}  =  x z_{xx} +(\alpha +1-x) z_x + b(x)u \
\ t \geq 0, \ \  x \in \RR_{+}, 
\end{equation}
where $b \in L^{2}(\RR_{+},\mu_{\alpha})$ and
the control $u$ belong to $L^{2}(0,t_1;\RR_{+})$. 

The equation (\ref{C111}) is approximately controllable if and only if  
$$ \int_{\RR_{+}} b(x) l^{\alpha}_{\nu}(x)x^{-\alpha}e^{x} dx \not= 0, \ \ \nu = 0,1,2,\dots. $$
In particular, if $\alpha=\frac{n}{2}-1$ then the equation
(\ref{C111}) is associated to the Cox-Ingersoll-Ross (CIR)
processes with a single control and therefore the controlled CIR
can never be exactly controllable on $[0,t_1]$.

\item[b)] The Jacobi equation in one variable  with a single control
\begin{equation}
\label{C222}
z_{t}  =  (1-x^2) z_{xx} +((\beta -\alpha-\left(\alpha +\beta +2\right)x) z_x + b(x)u \
\ t \geq 0, \ \  x \in [-1,1], 
\end{equation}
where $b \in L^{2}([-1,1],\mu_{\alpha,\beta})$ and
the control $u$ belong to $L^{2}(0,t_1;[-1,1])$. 

The equation (\ref{C222}) is approximately controllable if and only if  
\begin{eqnarray*}
 \int_{[-1,1]} b(x)p^{\alpha,\beta}_{\nu} \left( 1-x\right) ^{-\alpha }\left( 1+x\right) ^{-\beta
} dx \not= 0, \, \nu = 0, 1,2,\ldots.
\end{eqnarray*}

\end{enumerate}

\end{eje}

\begin{obs}
\label{obs2}
Notice that in each case,  the approximated controllability is totally determined by the non-orthogonality of the function $b \in L^{2}(\RR_{+},\mu_{\alpha})$ (respectively, $b\in L^{2}([-1,1],\mu_{\alpha,\beta})$) and the Laguerre (respectively, Jacobi) polynomials and it is independent of choice of control $u$.
\end{obs}

Finally, we will make some comments about the controllability  of general Sturm-Liouville equations. 
From a general point of view our arguments require  of the following ingredients:
\begin{enumerate}
\item A measure space $(\Omega, \Sigma, \mu)$, where $\Omega\subseteq \CC^{d}$ and $\mu$ is a Borel measure defined on $\Omega$.
\item An differential operator Sturm-Liouville type $\mathcal{L}$, whose eigenfunctions $\{y_{n}\}_{n\geq 0}$ form a complete orthogonal system  in $L^{2}(\Omega, d\mu)$ with complex eigenvalues  $\{\lambda_{n}\}_{n\geq 0}$ such that $\Re(\lambda_{n})\rightarrow \infty$ as $n \rightarrow \infty$.
\item A sequence of orthogonal projections $\{P_{n}\}_{n\geq 0}$ associated to the complete orthogonal system $\{y_{n}\}_{n\geq 0}$.
\item The  Hilbert space of  complex square sumable sequences $l^{2}$.
\end{enumerate}

With these ingredients the semigroup of operators $\{T_{t}\}_{t\geq}$ given by
$$T_{t}f=\sum_{n\geq 0}e^{-\lambda_{n}t}P_{n}f$$
is a strongly continuous  semigroup of compact operators,having infinitesimal generator,
$$\mathcal{L} =  \sum_{n\geq 0}(-\lambda_{n})P_{n}f,$$
with domain
$$D(\mathcal{L})=\left\{f\in
L^{2}(\Omega, d\mu): \sum_{n\geq 0} \|\lambda_{n}P_n f\|_{L^{2}(\Omega, d\mu)}^{2}< \infty\right\}.$$
Then for $b\in L^{2}(\Omega, d\mu)$ fixed, we consider the linear and bounded operator \\$B:l^{2}\rightarrow L^{2}(\Omega, d\mu)$ defined by
$$BU=\sum_{n\geq 0}U_{n}\langle b,y_{n}\rangle_{L^{2}(\Omega, d\mu)}y_{n}.$$

Then, the controlled equation associated to the Sturm-Liouville differential operator $\mathcal{L}$, 
$$z^{\prime}(t)=\mathcal{L}z(t)+Bu(t), \quad t\geq 0$$
is approximately controllable on $[0,t_1]$, if and only if,
$$P_{n}b\not=0, \quad \mbox{ for all } n\geq 0.$$

The special case $d=1$ and $\Omega$ be the unit circle of the complex plane, the support of the orthogonality measure $\mu$ for the Besell polynomials $\{B_{n}\}_{n\geq 0}$,  which are eigenfunctions of the differential operator
\begin{eqnarray}
\mathcal{L}&=& x^{2}\frac{d^{2}}{dx^{2}}+(2x+2)\frac{d}{dx},
\end{eqnarray}
with eigenvalue $n(n+1)$, $n\geq 0$.  Since these polynomials constitute a complete orthogonal system  in $L^{2}(\Omega, d\mu)$, then, if we consider the Bessel equation with a single control
\begin{equation}
\label{C223}
z_{t}  =  x^{2} z_{xx} +(2x+2)z_x + b(x)u \
\ t \geq 0, \ \  x \in \Omega, 
\end{equation}
where $b \in L^{2}(\Omega, d\mu)$ and
the control $u$ belong to $L^{2}(0,t_1;\Omega)$, we have that (\ref{C223}) is approximately controllable if and only if  
\begin{eqnarray*}
\langle b,B_{n}\rangle_{L^{2}(\Omega, d\mu)} \not= 0, \quad \mbox{ for all } n\geq 0.
\end{eqnarray*}

\vspace{2cm}

\begin{tabular}{cc}
\small
\parbox{7cm}{
Diomedes B\'arcenas, \, Hugo Leiva\\
Departamento de Matem\'aticas\\
Universidad de Los Andes\\
M\'erida 5101\\
VENEZUELA\\
\quad\\}
& \small
\parbox{7cm}{
\quad\\
Yamilet Quintana\\
Departamento de Matem\'aticas\\
Apartado Postal: 89000, Caracas 1080 A\\
Universidad Sim\'on Bol\'{\i}var\\
VENEZUELA\\
\quad\\
\quad\\} 
\end{tabular}

\vspace{-.5cm}

\begin{tabular}{cc}
\small
\parbox{7cm}{
Wilfredo Urbina\\
Departamento de Matem\'aticas
Facultad de Ciencias\\
Universidad Central de Venezuela, \\
Caracas, VENEZUELA\\
\quad\\}
& \small
\parbox{7cm}{
\quad\\
and Department of Mathematics and Statistics,\\
University of New Mexico, Albuquerque, New Mexico, 8713, USA\\
\quad\\
\quad\\}

\end{tabular}
\end{document}